\documentclass{amsart}
\usepackage{amssymb, amsfonts}
\usepackage[v2]{xy}

\newtheorem{thm}{Theorem}[section]
\newtheorem{cor}[thm]{Corollary}
\newtheorem{prop}[thm]{Proposition}
\newtheorem{lem}[thm]{Lemma}

\theoremstyle{definition}
\newtheorem{defn}[thm]{Definition}

\newtheorem{exmp}[thm]{Example}

\theoremstyle{remark}
\newtheorem{rem}[thm]{Remark}
\newtheorem{rems}[thm]{Remarks}

\DeclareFontFamily{OMS}{rsfs}{\skewchar\font'60}
\DeclareFontShape{OMS}{rsfs}{m}{n}{<-5>rsfs5 <5-7>rsfs7 <7->rsfs10 }{}
\DeclareSymbolFont{rsfs}{OMS}{rsfs}{m}{n}
\DeclareSymbolFontAlphabet{\scr}{rsfs}

\newcommand{\bF}{\mathbb{F}}

\newcommand{\bZ}{\mathbb{Z}}

\newcommand{\al}{\alpha}
\newcommand{\be}{\beta}

\newcommand{\et}{\eta}
\newcommand{\io}{\iota}

\newcommand{\si}{\sigma}

\newcommand{\SI}{\Sigma}

\newcommand{\OM}{\Omega}

\newcommand{\com}{\circ}     
\newcommand{\iso}{\cong}     
\newcommand{\htp}{\simeq}    
\newcommand{\ten}{\otimes}   
\newcommand{\sma}{\wedge}    

\newcommand{\rtarr}{\longrightarrow}

\def\quickop#1{\expandafter\newcommand\csname #1\endcsname{\operatorname{#1}}}
\quickop{Hom} \quickop{End} \quickop{Aut} \quickop{Tel} \quickop{Mic} 
\quickop{Ext} \quickop{Tor} \quickop{Id} \quickop{Coker} \quickop{Ker}
\quickop{Lim} \quickop{Colim} \quickop{Holim} \quickop{Hocolim}
\quickop{id} \quickop{tel} \quickop{mic} \quickop{coker} \quickop{Map}
\quickop{colim} \quickop{holim} \quickop{hocolim} \quickop{im} 

\makeatletter
\let\c@equation\c@thm
\makeatother
\numberwithin{equation}{section}

\def\:{\colon}
\def\<{\left\langle}
\def\>{\right\rangle}

\def\C{\mathbb{C}}

\def\F{\mathbb{F}}
\def\H{\mathbb{H}}

\def\R{\mathbb{R}}
\def\Z{\mathbb{Z}}
\def\BoP{{BoP}}
\def\BP{BP}
\def\MU{MU}

\def\eo{eo}
\def\ko{ko}
\def\ku{ku}
\def\SU{\mathrm{SU}}
\def\MSU{M\SU}
\def\CP{\C\mathrm{P}}
\def\HP{\H\mathrm{P}}
\def\RP{\R\mathrm{P}}

\DeclareMathOperator{\indeter}{indeter}

\title{Minimal atomic complexes}

\author{A.J. Baker and J.P. May}

\address{Mathematics Department, University of Glasgow,
Glasgow G12 8QW, Scotland.}
\email{a.baker@maths.gla.ac.uk}

\address{
Department of Mathematics, University of Chicago,
Chicago, IL 60637, USA.}
\email{may@math.uchicago.edu}



\begin{document}

\begin{abstract}
Hu, Kriz and May recently reexamined ideas implicit in Priddy's elegant
homotopy theoretic construction of the Brown-Peterson spectrum at a prime~$p$.
They discussed May's notions of \emph{nuclear complexes} and of \emph{cores} of
spaces, spectra, and commutative $S$-algebras. Their most striking conclusions,
due to Hu and Kriz, were negative: cores are not unique up to equivalence, and 
$\BP$ is not a core of $\MU$ considered as a commutative $S$-algebra, although it 
is a core of $\MU$ considered as a $p$-local spectrum. We investigate these ideas 
further, obtaining much more positive conclusions. We show that nuclear complexes 
have several non-obviously equivalent characterizations. Up to equivalence, they
are precisely the {\em irreducible complexes}, the {\em minimal atomic complexes}, 
and the Hurewicz complexes with trivial mod $p$ Hurewicz homomorphism above the 
Hurewicz dimension, which we call {\em complexes with no mod $p$ detectable
homotopy}. Unlike the notion of a nuclear complex, these other notions are all 
invariant under equivalence. This simple and conceptual criterion for a complex
to be minimal atomic allows us to prove that many familiar spectra, such as $\ko$, 
$\eo_2$, and $\BoP$ at the prime $2$,  all $\BP\<n\>$ at any prime $p$, and the 
indecomposable wedge summands of $\Sigma^\infty\CP^\infty$ and $\Sigma^\infty\HP^\infty$ 
at any prime $p$ are minimal atomic. 
\end{abstract}

\maketitle

\section*{Introduction}

Atomic spaces and spectra have long been studied. They are so tightly bound together
that a self-map which induces a isomorphism on homotopy in the Hurewicz dimension
must be an equivalence. Atomic spaces and spectra can often be shrunk to ones with
smaller homotopy groups. Minimal ones can be shrunk no further. Clearly, these are
very natural objects of study. They seem to have been
first introduced in \cite{HKM}. Spheres, $2$-cell complexes that are not wedges, and 
$K(\pi,n)$'s for cyclic groups $\pi$ are obviously minimal atomic, but there are many 
much more interesting examples.

Nuclear complexes are atomic complexes (spaces or spectra) that are built up in an 
especially economical way. They are minimal atomic, and we shall see that every minimal
atomic complex is equivalent to a nuclear complex. We regard the invariant notion of a minimal 
atomic complex as the more fundamental one. The combinatorial notion of a nuclear complex 
provides us with a tool for proving things about minimal atomic complexes.  

We think of atomic complexes as analogues of ``atomic modules'', namely modules for which 
a non-trivial self-map is an isomorphism. We think of minimal atomic complexes as analogues 
of irreducible modules.  We give a definition of an irreducible complex that makes 
this analogy more transparent, and we prove that the irreducible complexes are precisely 
the minimal atomic complexes. In one direction, the implication gives a homotopical 
analogue of Schur's lemma. Just as in algebra, we suggest that the irreducible, or 
minimal atomic, complexes are more basic mathematical objects than the atomic complexes. 
Carrying the analogy further, we see that the spectra we consider admit (dual) 
composition series of infinite length, constructed in terms of irreducible complexes.

There is a different and much more elementary notion of minimality, implicitly due to 
Cooke \cite{Cooke}, such that {\em any} complex is equivalent to a minimal complex. 
This notion is also combinatorial and noninvariant. We prove that a Hurewicz complex 
that is minimal in this sense is nuclear if and only if it has no mod $p$ detectable 
homotopy. It is well known that the latter condition implies wedge indecomposability 
(e.g. \cite[5.4]{Pengelley}). The fact that it implies irreducibility is much 
stronger.  We also prove a converse, leading to the conclusion that a complex 
is irreducible, or equivalently minimal atomic, if and only if it has no mod $p$\, detectable 
homotopy. This allows us to show that a variety of familiar spectra are in fact 
minimal atomic. 

\section{Definitions and statements of results}

Here we give the precise definitions needed to make sense of the introduction
and state our main theorems.  We write things so that the stable reader can view these as 
statements about spectra, and the unstable reader can view them as 
statements about (based) spaces. We adopt the following conventions throughout. 
They allow us to treat 
spaces and spectra uniformly and to avoid repeated mention of the fact that we are working 
$p$-locally under connectivity and finite type hypotheses. 

We agree once and for all that all spaces and spectra $X$ are to be localized at a fixed prime $p$. 
Thus $S^n$, for example, means a $p$-local sphere. We also agree that all spaces and spectra are 
to be $p$-local CW spaces or spectra, so that the domains of their attaching maps are $p$-local spheres. 
Spaces are to be simply connected, and their attaching maps are to be based. Spectra are to be bounded 
below.  In either case, we say that $X$ {\em has Hurewicz dimension $n_0$} if $X$ is $(n_0-1)$-connected, 
but not $n_0$-connected. Thus $n_0\geqslant 2$ in the case of spaces, and there is no real loss of generality 
if we take $n_0=0$ in the case of spectra.  We may assume without loss of generality that there are no 
cells (except the base vertex) of dimension less than $n_0$.  We assume further that there are only 
finitely many cells in each dimension. We agree to use the ambiguous term ``complex'' to mean such a 
$p$-local CW space or spectrum.  We say that $X$ is a {\em Hurewicz complex} if it has a single cell
in dimension $n_0$. We write $X_n$ for the $n$-skeleton of a complex $X$. Thus $X_{n+1}$ is the cofiber 
of a map $j_n: J_n\rtarr X_n$, where $J_n$ is a finite wedge of ($p$-local) $n$-spheres $S^n$. If $X$
is a Hurewicz complex, $X_{n_0} = S^{n_0}$. We shall use these notations generically. 

By $H_*(X)$, we always understand (reduced) homology with $p$-local coefficients. Any 
$(n_0-1)$-connected space or spectrum such that each $H_{n}(X)$ is a finitely generated 
$\bZ_{(p)}$-module is weakly equivalent to a complex in the sense that we have just specified.
If, further, $H_{n_0}(X;\bF_p)=\bF_p$ or, equivalently, $\pi_{n_0}(X)$ is a cyclic 
$\bZ_{(p)}$-module, then $X$ is weakly equivalent to a Hurewicz complex.

We begin with definitions of concepts that are invariant under equivalence
and the statement of our main characterization theorem relating them.

\begin{defn}\label{defn:inv} Consider complexes $X$ and $Y$ of Hurewicz dimension $n_0$. 
Think of $Y$ as fixed but $X$ as variable.
\begin{enumerate}
\item[(i)] A map $f: X\rtarr Y$ is a {\em monomorphism} if 
$f_*\: \pi_{n_0}(X)\ten \bF_p \rtarr \pi_{n_0}(Y)\ten\bF_p$
and all $f_*\: \pi_n(X)\rtarr \pi_n(Y)$ are monomorphisms
\item[(ii)] $Y$ is {\em irreducible} if any monomorphism $f: X\rtarr Y$
is an equivalence.
\item[(iii)] $Y$ is {\em atomic} if it is a Hurewicz complex and a self-map 
$f:Y\rtarr Y$ that induces an isomorphism on $\pi_{n_0}$ is an equivalence. 
\item[(iv)] $Y$ is {\em minimal atomic} if it is atomic and any monomorphism 
$f: X\rtarr Y$ from an atomic complex $X$ to $Y$ is an equivalence.
\item[(v)] $Y$ {\em has no mod $p$ detectable homotopy} if $Y$ is a Hurewicz
complex and the mod $p$ Hurewicz homomorphism $h: \pi_n(Y)\rtarr H_n(Y;\bF_p)$ 
is zero for all $n>n_0$. 
\item[(vi)] $Y$ is {\em $H^*$-monogenic}\, if $H^*(Y;\bF_p)$ is a cyclic algebra
(in the case of spaces) or module (in the case of spectra) over the mod $p$\,
Steenrod algebra $A$.
\end{enumerate}
\end{defn}

\begin{rems}\label{rems:many} 
We offer several comments on these notions.
\begin{enumerate}
\item[(i)]  The structure theory for finitely generated modules over 
a PID implies that if $f\:X \rtarr Y$ is a monomorphism, then 
$f_*(\pi_{n_0}(X))$ is a direct summand of $\pi_{n_0}(Y)$. If $X$ is 
a Hurewicz complex, this summand is cyclic. If $X$ and $Y$ are both 
Hurewicz complexes, then $f$ induces an isomorphism on $\pi_{n_0}$.
\item[(ii)]  In \cite[1.1]{HKM}, following \cite{Wilk} and other early 
sources, $Y$ was defined to be irreducible if it has no non-trivial 
retracts. On the space level, that concept has its uses, but we think 
that ``irreducible'' is the wrong name for it. We suggest 
``irretractible''. On the spectrum level, irretractibility is equivalent 
to wedge indecomposability. However, just as in algebra, irreducibility 
should be stronger rather than weaker than atomic. That is, there should be 
implications irreducible $\Longrightarrow$ atomic $\Longrightarrow$ 
indecomposable. One could avoid the conflict with the earlier literature 
by using the word ``simple'' instead of ``irreducible'', the two being 
synonymous in algebra, but that risks confusion with the standard use of 
the term ``simple'' in topology.
\item[(iii)] A complex can well have more than one cell in its Hurewicz 
dimension and still have the property that a self-map that induces an 
isomorphism on $\pi_{n_0}$ is an equivalence. A particularly interesting 
example is given in \cite[\S4]{AdamsKuhn}. It might be sensible to delete 
the requirement that $Y$ be a Hurewicz complex from the definition of atomic. 
By Theorem \ref{thm:core} below, the notion of minimal atomic would not change.
\item[(iv)] Since our methods are cellular, we definitely mean to consider 
$p$-local rather than $p$-complete spaces and spectra. However, 
Definition \ref{defn:inv} makes just as much sense in the $p$-complete case as 
the $p$-local case, and it is well worth studying there. Since a finite type 
$p$-complete space or spectrum is the $p$-completion of a finite type $p$-local 
space or spectrum, one can easily deduce conclusions in the $p$-complete case 
from the results here. We leave the details to the interested reader.
\item[(v)] A Hurewicz complex $Y$ has no mod $p$\, detectable homotopy if and only if
there are no permanent cycles in dimension greater than $n_0$ on the zeroth 
row of the classical (unstable or stable) mod $p$\, Adams spectral sequence 
for $Y$. It is a much more computable condition than the others.
\end{enumerate}
\end{rems} 

\begin{thm}[The characterization theorem]\label{thm:main} 
The following conditions on a complex $Y$ are equivalent.
\begin{enumerate}
\item[(i)] $Y$ is irreducible.
\item[(ii)] $Y$ is minimal atomic.
\item[(iii)] $Y$ has no mod $p$ detectable homotopy.
\end{enumerate}
\end{thm}

The fact that irreducible complexes are atomic should be viewed as a homotopical 
analogue of Schur's lemma since its intuitive content is that a non-trivial self-map of an
irreducible complex must be an equivalence. Of course, it is consistent with the analogy
that not all atomic complexes are irreducible. When \cite{HKM} was written, examples of 
minimal atomic spectra seemed hard to come by. The following result now gives us many 
interesting examples.

\begin{cor}\label{cor:monat}
If $Y$ is $H^*$-monogenic, then $Y$ has no mod~$p$\, detectable homotopy
and is therefore minimal atomic.
\end{cor}
\begin{proof}
For spectra, the long exact sequence on $\Ext$ arising from the 
epimorphism $\SI^{n_0} A\rtarr H^*(Y;\F_p)$ implies that
the zeroth row $\Hom_{A}^{0,*}(H^*(Y;\F_p),\F_p)$ of the Adams spectral 
sequence is $\bF_p$ concentrated in degree $n_0$. Similarly for spaces.
\end{proof}

\begin{rem} The converse of Corollary \ref{cor:monat} fails. For $q\geqslant2$, 
Moore spaces and spectra $M(\bZ/p^q,n)$ and Eilenberg-Mac\,Lane spaces and 
spectra $K(\bZ/p^q,n)$ are elementary examples of complexes that are minimal 
atomic but not $H^*$-monogenic.
\end{rem}

The proof of Theorem \ref{thm:main} proceeds by
first showing in Theorem \ref{thm:first} that (i) and (ii) are each equivalent to a statement about 
one noninvariant cellular construction and then showing in Theorem \ref{thm:second} that (ii)
and (iii) are each equivalent to a statement about another noninvariant cellular construction. 
The first noninvariant construction, codified in Theorems \ref{thm:core} 
and \ref{thm:Min}, is based on the notions of nuclear complexes and cores 
introduced by Priddy \cite{Priddy-BP} and May \cite{HKM}.

\begin{defn} A {\em nuclear complex} is a Hurewicz complex $X$ such that
\begin{equation}\label{key}
\Ker(j_{n*}: \pi_n(J_n) \rtarr \pi_n(X_n))\subset p\cdot \pi_n(J_n)
\end{equation}
for each $n$. Observe that $X$ is nuclear if and only if each $X_n$ for $n\geqslant n_0$
is nuclear. A {\em core} of a complex $Y$ is a nuclear complex $X$ together with a monomorphism 
$f: X\rtarr Y$.
\end{defn}

This notion of a core is more general than that of \cite[1.7]{HKM}, 
where it was assumed that $\pi_{n_0}(Y)$ is cyclic; that is, the definition
there was restricted to Hurewicz complexes $Y$. Since cores are not unique
even when $Y$ is a Hurewicz complex, the more general notion seems 
preferable. Of course, with our perspective that cores are analogues 
of irreducible sub-modules, the non-uniqueness is only to be expected.
With the present language, the following results are proven 
in \cite[1.5, 1.6]{HKM}.

\begin{prop}\label{prop:partial} A nuclear complex is atomic.
\end{prop}

\begin{thm}\label{thm:core} 
If $Y$ has Hurewicz dimension $n_0$ and $C$ is a cyclic direct
summand of $\pi_{n_0}(Y)$, then there is a core $f:X\rtarr Y$ such that
$f_*(\pi_{n_0}(X))=C$.
\end{thm}

This allows us to construct homotopical analogues of composition series. That is, 
for spectra, we can shrink homotopy groups inductively by successively taking 
cofibers of cores, as discussed in \cite[1.8]{HKM}. This works less well for spaces, 
where we would have to take fibers and so gradually decrease the Hurewicz dimension.
To make the analogy with algebra precise, recall that a (countably infinite) composition 
series of a module $Y$ is a sequence of monomorphisms
$$\xymatrix@1{
Y = Y_0 & Y_1 \ar[l]_-{i_0} & \cdots \ar[l] & Y_n \ar[l] & Y_{n+1} \ar[l]_{i_n} & \ar[l] \cdots\\}$$
such that the cokernels of the $i_n$ are irreducible and $\lim Y_n = 0$. A {\em dual composition series}\,
of $Y$ is a sequence of epimorphisms
$$\xymatrix@1{
Y = Y_0 \ar[r]^-{p_0} & Y_1 \ar[r]  & \cdots  \ar[r] & Y_n \ar[r]^-{p_n} & Y_{n+1} \ar[r] & \cdots\\}$$
such that the kernels of the $p_n$ are irreducible and $\colim Y_n = 0$. 
For a spectrum $Y$, we construct an analogous sequence by letting 
$p_n: Y_n\rtarr Y_{n+1}$ be the cofiber of a core $f_n: X_n\rtarr Y_n$. 
Each $p_n$ induces an epimorphism on all homotopy groups, we kill 
$\pi_{n_0}(Y)$ in finitely many steps, then kill $\pi_{n_0+1}(Y)$ in 
finitely many steps, and so on. If\, $Y$ has non-zero homotopy groups 
in only finitely many dimensions, then this sequence has only finitely 
many terms.

Theorem \ref{thm:core} has the following immediate consequence.

\begin{cor}\label{cor:oneway}  A core of a minimal atomic complex is an 
equivalence, hence a minimal atomic complex is equivalent to a nuclear 
complex.
\end{cor}

The converse, which completes Proposition \ref{prop:partial}, was conjectured 
in \cite[1.12]{HKM}. We will prove it in Section \ref{sec:Min}.

\begin{thm}\label{thm:Min}
A nuclear complex is a minimal atomic complex.
\end{thm}

\begin{thm}\label{thm:first} Conditions \emph{(i)} and \emph{(ii)} of \emph{Theorem \ref{thm:main}} 
are each equivalent to the condition that any core $f: X \rtarr Y$ of\, $Y$ is an equivalence.
\end{thm}
\begin{proof} Corollary \ref{cor:oneway} and Theorem \ref{thm:Min} show that $Y$
is minimal atomic if and only if any core of $Y$ is an equivalence. If $Y$ is 
irreducible and $f:X\rtarr Y$ is a core, then $f$ is an equivalence by the
definition of irreducibility. Conversely, suppose that $Y$ is minimal atomic and
let $f:X\rtarr Y$ be a monomorphism. Let $g:W\rtarr X$ be a core of $X$. Then the 
composite $f\com g:W\rtarr Y$ is a core of $Y$ and therefore an equivalence. This implies
that $f$ induces an epimorphism and hence an isomorphism on homotopy groups. Thus
$f$ is an equivalence.
\end{proof}

\begin{rem} With these implications in place, it is perhaps better to redefine
the notion of core invariantly, taking $X$ to be minimal atomic but not 
necessarily nuclear. There is no substantive difference.
\end{rem}

To tie in the Hurewicz homomorphism condition (iii) of Theorem \ref{thm:main}, 
we use another, very different, noninvariant notion of minimality for a complex $X$. 
Of course, our complexes have $p$-local chain complexes specified by $C_n(X) = H_n(X_{n}/X_{n-1})$.

\begin{defn} A complex $X$ is {\em minimal} if the differential on its mod $p$ chain 
complex $C_*(X;\bF_p)$ is zero. It is {\em minimal Hurewicz} if it is minimal and 
Hurewicz. Observe that $X$ is minimal if and only if each $X_n$ is minimal.
\end{defn}

A simple inductive argument gives a homological reformulation of this notion.

\begin{lem}\label{lem:cute} 
A complex $X$ is minimal if and only if the inclusion of skeleta 
$X_n \rtarr X_{n+1}$ induces an isomorphism
$$ H_n(X_n;\F_p)\rtarr H_n(X_{n+1};\F_p)=H_n(X;\F_p)$$
for each $n$.
\end{lem}

The following two results codify our second noninvariant construction. The first is 
implicit in Cooke's paper \cite[Theorem A]{Cooke}, which gives an integral space level 
version. Cooke described the result as ``a well-known, basic fact''. For a recent 
reappearance, see \cite[4.C.1]{Hatcher}. The proof is very easy, but we shall run 
through it in Section \ref{sec:Cooke} in view of the importance of the result to our 
work.

\begin{thm}\label{thm:Cooke}
For any complex $Y$, there is a minimal complex $X$ and an
equivalence $f: X\rtarr Y$.
\end{thm}

We prove the following result in Section \ref{sec:Hur}. 

\begin{thm}\label{thm:Hur} If $X$ is a nuclear complex, then $X$ 
has no mod $p$ detectable homotopy. If $X$ is a minimal Hurewicz 
complex, then $X$ is nuclear if and only if it has no mod $p$ 
detectable homotopy.
\end{thm}

\begin{thm}\label{thm:second} Conditions \emph{(ii)} and \emph{(iii)} 
of \emph{Theorem \ref{thm:main}} are each 
equivalent to the condition that any equivalence $X\rtarr Y$ from a minimal complex $X$ to $Y$
is a core of $Y$; that is, a minimal complex equivalent to $Y$ is nuclear.
\end{thm}
\begin{proof} Since a minimal atomic complex is equivalent to a nuclear complex, the 
second statement of Theorem \ref{thm:Hur} implies that (iii) is equivalent to the condition 
specified in the statement and that (iii) implies (ii).  Similarly, the first statement of 
Theorem \ref{thm:Hur} implies that (ii) implies (iii).  
\end{proof}

Corollary \ref{cor:monat} and, more generally, the implication (iii) implies (ii) 
provide a powerful tool for detecting minimal atomic complexes. We give some 
general results that illustrate the use of this criterion in the next section.

Restricting to spectra, we see in Section~\ref{sec:SomeNuclSpectra} that $\ko$ and 
$\eo_2$ at the prime $2$, $\BP\langle n\rangle$ at any prime $p$, and the indecomposable 
wedge summands of
$\Sigma^\infty\CP^\infty$ and $\Sigma^\infty\HP^\infty$ at any prime $p$ are minimal 
atomic. We give a few other examples and remarks, but we regard this section as just 
a beginning. Our results imply that minimal atomic complexes exist in abundance, and
something closer to a classification of them would be desirable.

In Section~\ref{sec:BoP}, we describe Pengelley's $2$-local spectrum $\BoP$ as a 
nuclear complex and thereby give it a new construction that is independent of
\cite{Pengelley}. This is in the same spirit as Priddy's construction of 
$\BP$ \cite{Priddy-BP}, which we recall in Section \ref{sec:SomeNuclSpectra}.  
The key step in the proof is deferred to Section~\ref{sec:\ko}.
A brief Appendix corrects minor errors in one of the proofs in ~\cite{HKM}.

\section{Constructions on minimal atomic complexes}

We indicate briefly how the collection of minimal atomic complexes 
behaves with respect to some basic topological constructions.
The proofs are direct consequences of the ``no mod $p$\,
detectable homotopy'' characterization of minimal atomic. 
The following triviality may help the reader see the
various implications.

\begin{lem}\label{lem:trivial} Consider a commutative diagram
$$\xymatrix{
A \ar[d]_h \ar[r]^-{f} & A' \ar[d]^{h'}\\
B  \ar[r]^-{g} & B'\\}$$
of Abelian groups. If $f$ is an epimorphism and $h=0$, then $h'=0$. If $g$ is a monomorphism and $h'=0$, then $h=0$.
\end{lem}

We begin by recording an immediate consequence of Theorems \ref{thm:Hur} and \ref{thm:second}.

\begin{prop} If $Y$ is minimal atomic, then $Y$ is equivalent to a 
complex $X$ whose skeleta $X_n$ for $n\geqslant n_0$ are minimal 
atomic.
\end{prop}

There is no reason to believe that the skeleta of $Y$ itself are
minimal atomic. We have a more invariant analogue for Postnikov 
sections, which we denote by $Y[n]$.

\begin{prop} A complex $Y$ is minimal atomic if and only if $Y[n]$ 
is minimal atomic for each $n\geqslant n_0$. 
\end{prop}
\begin{proof}  We have the following commutative diagram.
$$\xymatrix{
\pi_q(Y) \ar[r] \ar[d]_{h} & \pi_q(Y[n]) \ar[d]^{h} \\
H_q(Y;\F_p) \ar[r]  & H_q(Y[n];\F_p)
\\}$$
Since $\pi_q(Y[n])$ is zero for $q>n$ and the horizontal
arrows are isomorphisms for $q\leqslant n$, the conclusion
is immediate from Lemma \ref{lem:trivial}.
\end{proof}

In the following result, which is only of interest for spaces, we 
consider the loop and suspension functors.

\begin{prop} If either $\OM Y$ or $\SI Y$ is minimal atomic, then 
so is $Y$.
\end{prop}
\begin{proof} This is immediate from the following commutative diagrams.
$$\xymatrix{
\pi_q(\OM Y) \ar[r]^-{\iso} \ar[d]_{h} & \pi_{q+1}(Y) \ar[d]^{h} \\
H_q(\OM Y;\F_p) \ar[r]_-{\si}  & H_{q+1}(Y;\F_p)\\}
\ \ \ \ \ \ \ 
\xymatrix{
\pi_q(Y) \ar[r]^-{\SI} \ar[d]_{h} & \pi_{q+1}(\SI Y) \ar[d]^{h} \\
H_q(Y;\F_p) \ar[r]_-{\iso}  & H_{q+1}(\SI Y;\F_p) \qed\\}$$
\renewcommand{\qed}{}
\end{proof}

This result has an analogue that relates minimal atomicity for spaces and spectra.
Here, exceptionally, we must distinguish the two contexts notationally.

\begin{prop} If $E$ is a spectrum of Hurewicz dimension $n_0\geqslant 2$ whose 
$0$th space $\OM^{\infty}E$ is minimal atomic, then $E$ is minimal atomic. 
If $Y$ is a simply connected space whose suspension spectrum $\SI^{\infty}Y$ 
is minimal atomic, then $Y$ is minimal atomic.
\end{prop}
\begin{proof} This is immediate from the following commutative diagrams.
$$\xymatrix{
\pi_q(\OM^{\infty} E) \ar[r]^-{\iso} \ar[d]_{h} & \pi_{q}(E) \ar[d]^{h} \\
H_q(\OM^{\infty}E ;\F_p) \ar[r]_-{\si}  & H_{q}(E;\F_p)\\}
\ \ \ \ \ \ \ 
\xymatrix{
\pi_q(Y) \ar[r] \ar[d]_{h} & \pi_{q}(\SI^{\infty} Y) \ar[d]^{h} \\
H_q(Y;\F_p) \ar[r]_-{\iso}  & H_{q}(\SI^{\infty} Y;\F_p) \qed\\}$$
\renewcommand{\qed}{}
\end{proof}

\section{The proof of Theorem \ref{thm:Min}}\label{sec:Min}

The proof in \cite[1.5]{HKM} that a nuclear complex $X$ is atomic starts with a self map
$f: X\rtarr X$ which is an isomorphism on $\pi_{n_0}(X)$ and deduces that $f$ is an
equivalence.  The cited proof readily adapts to give the following analogue.

\begin{prop}\label{prop:Nuclear-CoreUnique}
Let $X$ and $Y$ be nuclear complexes of Hurewicz dimension $n_0$ and let $f\:X\rtarr Y$ 
be a core of\, $Y$. Then~$f$ is an equivalence.
\end{prop}
\begin{proof}
Take $f\:X\rtarr Y$ to be cellular. Since $f$ is a monomorphism between Hurewicz 
complexes, $f:X_{n_0}\rtarr Y_{n_0}$ is an equivalence. Assume that $f:X_n\rtarr Y_n$ 
is an equivalence. We must show that $f:X_{n+1}\rtarr Y_{n+1}$ is an equivalence. The 
attaching maps of $X$ and $Y$ give rise to the following map of cofibre sequences.
$$\xymatrix{
J_n \ar[r]^-{j_n} \ar[d]_{f} & X_n \ar[r] \ar[d]^{f} & X_{n+1} \ar[d]^{f}\\
K_n\ar[r]_-{k_n} & Y_n \ar[r] & Y_{n+1}\\}$$
Passing to homology, this gives rise to a commutative diagram with exact rows.
$$\xymatrix{
0 \ar[r] & H_{n+1}(X_{n+1}) \ar[r] \ar[d]_{f_*} & H_n(J_n) \ar[r]^-{{j_n}_*} \ar[d]^{f_*}  
& H_n(X_n) \ar[r] \ar[d]^{f_*}_{\iso} & H_n(X_{n+1})\ar[r] \ar[d]^{f_*} & 0 \\
0 \ar[r] & H_{n+1}(Y_{n+1}) \ar[r] & H_n(K_n) \ar[r]_{{k_n}_*} & H_n(Y_n)\ar[r] 
& H_n(Y_{n+1}) \ar[r] & 0
\\}$$
It suffices to prove that the left and right vertical arrows are isomorphisms. 
By the five lemma and the Hurewicz theorem, this holds if $f_*\:\pi_n(J_n)\rtarr\pi_n(K_n)$ 
is an isomorphism. To see that this is so, consider the following diagram.
$$\xymatrix{
\pi_n(J_n) \ar[r]^-{{j_n}_*} \ar[d]_{f_*} & \pi_n(X_n) \ar[r] \ar[d]^{f_*} 
& \pi_n(X_{n+1})\ar[r] \ar[d]^{f_*} & 0 \\
\pi_n(K_n) \ar[r]_-{{k_n}_*} & \pi_n(Y_n) \ar[r] & \pi_n(Y_{n+1})\ar[r] & 0
\\}$$
The rows are exact, and a chase of the diagram shows that the right arrow $f_*$ is an epimorphism. 
Now consider the following diagram.
$$\xymatrix{
\pi_n(X_{n+1})\ar[d]_{f_*}  \ar[r]^-{\iso} & \pi_n(X) \ar[d]^{f_*}\\
\pi_n(Y_{n+1}) \ar[r]_-{\iso}  & \pi_n(Y)
\\}$$
Since its right arrow $f_*$ is a monomorphism, its left arrow $f_*$ is a monomorphism and
therefore an isomorphism. This implies that the right vertical arrow is an isomorphism in the 
following diagram.
$$\xymatrix{
0\ar[r] & \Ker\,{j_n}_* \ar[r]^-{i} \ar[d] & \pi_n(J_{n+1}) \ar[r]^-{\iso} \ar[d]^{f_*} 
& \text{Im}\, {j_n}_*\ar[r] \ar[d]^{\iso}  & 0\\
0\ar[r] & \Ker\,{k_n}_* \ar[r]_-{i} & \pi_n(K_{n+1}) \ar[r]_-{\iso} & \text{Im}\, {k_n}_*\ar[r] & 0
\\}$$
In view of (\ref{key}), both maps $i$ become $0$ after tensoring with $\bF_p$. This implies
that $f_*\ten \bF_p$ is an isomorphism, and therefore so is $f_*$. 
\end{proof}

\begin{proof}[Proof of Theorem \ref{thm:Min}]
Let $Y$ be a nuclear complex of Hurewicz dimension $n_0$ and let 
$f\:X\rtarr Y$ be a monomorphism, where $X$ is atomic. The same
argument as in the last part of the proof of Theorem \ref{thm:first}
shows that $f$ is an equivalence.
\end{proof}

\section{The proof of Theorem \ref{thm:Hur}}\label{sec:Hur}

We start with the following result, which is based on an observation 
of Priddy~\cite{Priddy-BP}. It gives a homological recasting of the definition 
of a nuclear complex.

\begin{lem}\label{lem:Observe} A Hurewicz complex of dimension $n_0$ is 
nuclear if and only if the mod $p$ Hurewicz homomorphism 
$h: \pi_{n}(X_{n})\rtarr H_{n}(X_n;\bF_p)$ is zero for $n > n_0$. 
\end{lem}
\begin{proof} Recall the defining property (\ref{key}) of a nuclear
complex. In the case of spaces, our assumption that $X$ is simply connected
allows us to quote the relative and absolute Hurewicz theorem to deduce that 
$$\pi_{n+1}(X_{n+1},X_n)\iso \pi_{n+1}(\SI J_n)\iso \pi_n(J_n)$$
from the trivial analogue in $p$-local homology.  In either the space or the 
spectrum context, we obtain the following commutative diagram with exact rows.
$$\xymatrix{
\pi_{n+1}(X_n) \ar[r] & \pi_{n+1}(X_{n+1}) \ar[r] \ar[d]^{h} 
& \pi_n(J_n) \ar[r]^-{j_*} \ar[d]^{h} & \pi_n(X_n) \ar[d]^{h}\\
0 \ar[r] & H_{n+1}(X_{n+1};\bF_p) \ar[r] & H_n(J_n;\bF_p) \ar[r]_-{j_*} & H_n(X_n;\bF_p)}$$
An easy diagram chase gives that (\ref{key}) holds for $n$ if and only if the left arrow $h$
is zero. The conclusion follows.
\end{proof}

We prove the following reformulation of Theorem \ref{thm:Hur} 
by relating this skeletal criterion to the Hurewicz homomorphism 
of $X$ itself.

\begin{prop}\label{prop:test}
Let $X$ be a Hurewicz complex of dimension $n_0$ and 
consider the following conditions.
\begin{enumerate}
\item[(i)] $X$ has no mod $p$ detectable homotopy.
\item[(ii)] $X$ is a minimal complex.
\end{enumerate}
If \emph{(i)} and \emph{(ii)} hold, then $X$ is nuclear. Conversely, 
if $X$ is nuclear, then \emph{(i)} holds.
\end{prop}
\begin{proof}
We have the following commutative diagram, where $n>n_0$.
$$\xymatrix{
\pi_n(X_n) \ar[r] \ar[d]_{h} & \pi_n(X) \ar[d]^{h} \\
H_n(X_n;\F_p) \ar[r]  & H_n(X;\F_p)
\\}$$
The conclusion is immediate from Lemmas \ref{lem:cute}, \ref{lem:trivial}, and \ref{lem:Observe}.
\end{proof}

\section{The proof of Theorem \ref{thm:Cooke}}\label{sec:Cooke}

We are given a complex $Y$. Recall that our complexes are simply connected in the 
case of spaces and bounded below in the case of spectra. More fundamentally, everything
is $p$-local. We have assumed that $H_*(Y)$ is of finite type, so that each $H_n(Y)$
is a direct sum of finitely many cyclic $\bZ_{(p)}$-modules $A_{n,i}$. We must construct 
a minimal complex $X$ and an equivalence $f:X\rtarr Y$, and it suffices for the latter to 
ensure that $f$ induces an isomorphism on $H_*$. The complex $X$ will have an $n$-cell 
$j_{n,i}$ for each free cyclic summand $A_{n,i}$ and an $n$-cell $j_{n,i}$ and an
$(n+1)$-cell $k_{n,i}$ with differential $q_i j_{n,i}$ for each summand $A_{n,i}$ of order
$q_i$. Since each $q_i$ must be a power of $p$, it will be immediate that the differential 
on $C_*(X;\bF_p)$ is zero.  The cells $j_{n,i}$ will map to cycles that represent the 
generators of the $A_{n,i}$, and the cells $k_{n,i}$ will map to chains with boundary
$q_if_*(j_{n,i})$.

Assume inductively that we have constructed the $n$-skeleton $X_n$ together with a 
(based) map $f_n: X_n\rtarr Y$ that induces an isomorphism on homology in dimensions 
less than $n$ and an epimorphism on $H_n$. More precisely, assume that $H_n(X_n)$ is
$\bZ_{(p)}$-free on basis elements given by cells $j_{n,i}$ that map to chosen generators of the
$A_{n,i}$. Let $Cf_n$ be the cofiber of $f_n$. Then $H_m(Cf)=0$ for $m\leqslant n$. The
kernel of $f_*: H_n(X_n)\rtarr H_n(Y)$ is free on the basis $q_i j_{n,i}$ for those
$i$ such that $A_{n,i}$ has finite order. These elements are the images of elements
$k''_{n,i}$ in $H_{n+1}(Cf_n)$, and $k''_{n,i} = h(k'_{n,i})$ for unique elements $k'_{n,i}$ 
in $\pi_{n+1}(Cf_n)$. Similarly, the chosen generators of the $A_{n+1,i}\subset H_{n+1}(Y)$
map to elements $j''_{n+1,i}\in H_{n+1}(Cf)$ with $j''_{n+1,i} = h(j'_{n+1,i})$. For spectra,
we have the connecting homomorphism $\pi_{n+1}(Cf_n) \rtarr \pi_n(X_n)$. For spaces, the 
relative Hurewicz theorem gives $\pi_{n+1}(Mf_n,X_n)\iso \pi_{n+1}(Cf)$, and we have
the connecting homomorphism $\pi_{n+1}(Mf_n,X_n) \rtarr \pi_n(X_n)$. Thus in either case
the elements $k'_{n,i}$ and $j'_{n+1,i}$ determine elements of $\pi_n(X_n)$. Choose maps
$S^n\rtarr X_n$ that represent these elements and use them as attaching maps for the 
construction of $X_{n+1}$ from $X_n$ by attaching cells $k_{n,i}$ and $j_{n+1,i}$.

Since the sequence $\pi_{n+1}(Cf_n)\rtarr \pi_n(X_n)\rtarr \pi_n(Y)$ is exact, 
these attaching maps become null homotopic in $Y$, and there is an extension 
$f_{n+1}: X_{n+1}\rtarr Y$ of $f_n$. Thus we can construct the following map of cofiber sequences.
$$\xymatrix{
X_n \ar @{=}[d] \ar[r]  &  X_{n+1} \ar[d]^{f_{n+1}} \ar[r] 
& X_{n+1}/X_n \ar[r] \ar[d] & \SI X_n \ar @{=} [d]\\
X_n \ar[r]_{f_n} & Y \ar[r] & Cf_n \ar[r] & \SI X_n \\}$$
This gives rise to the following commutative diagram with exact rows.

{\small
$$\xymatrix{
0 \ar[r]  &  H_{n+1}(X_{n+1}) \ar[d]_{(f_{n+1})_*} \ar[r] 
& H_{n+1}(X_{n+1}/X_n) \ar[r] \ar[d] & H_n(X_n) \ar @{=} [d] \ar[r] 
& H_n(X_{n+1}) \ar[d]^{(f_{n+1})_*} \ar[r] & 0\\
0 \ar[r] & H_{n+1}(Y) \ar[r] & H_{n+1}(Cf_n) \ar[r] & H_n(X_n) \ar[r] & H_n(Y) \ar[r] & 0\\}$$
}
Of course, the differential on $C_{n+1}(X_{n+1})$ is the composite
$$H_{n+1}(X_{n+1}/X_n) \rtarr H_n(X_n) \rtarr  H_n(X_n/X_{n-1}),$$
where the second arrow is a monomorphism. By construction, the first arrow sends 
the basis elements $k_{n,i}$ to $q_ij_{n,i}$ and the basis elements $j_{n+1,i}$ to 
zero, so that $H_{n+1}(X_{n+1})$ is $\bZ_{(p)}$-free on the basis elements $j_{n+1,i}$. 
By construction and a chase of the diagram, the map $f_{n+1}$ induces an isomorphism on 
$H_n$ and sends the basis elements $j_{n+1,i}$ to generators of the groups $A_{n+1,i}$. 
This completes the inductive step in the construction of $f:X\rtarr Y$.

\section{Spectrum level examples}\label{sec:SomeNuclSpectra}

As a first example, we revisit Priddy's construction 
\cite{Priddy-BP} of a nuclear spectrum equivalent to $\BP$. 
Although it was the motivating example for \cite{HKM}, it 
was not explicitly discussed there. We work with $p$-local
spectra in this section. Unless otherwise stated, $p$ is unrestricted.

\begin{exmp} $\BP$ is a minimal atomic spectrum, hence the 
canonical monomorphism $\BP\rtarr \MU$ is a core of $\MU$.
\end{exmp}
\begin{proof}
$H^*(\BP;\bF_p)$ is a cyclic $A$-module, hence Corollary \ref{cor:monat} applies.
\end{proof}

\begin{prop}\label{prop:BP}
Let $X$ be the nuclear complex of~\cite{Priddy-BP} defined by
starting with $S^0$ and inductively killing the homotopy groups in
odd degrees. Then there is an equivalence $X\rtarr\BP$.
\end{prop}
\begin{proof}
A minimal complex equivalent to $\BP$ has cells only in even degrees and is nuclear.
By construction, $X$ also has cells only in even degrees and is nuclear, and its 
non-zero homotopy groups only occur in even degrees. Obstruction theory gives maps 
$f\:X\rtarr\BP$ and  $g\:\BP\rtarr X$ that extend the identity on the bottom
cell. The composites $g\com f\:X\rtarr X$ and $f\com g\:\BP\rtarr\BP$ are equivalences 
since $X$ and $\BP$ are atomic.
\end{proof}

Recall that, for an odd prime $p$, there is a splitting of $\ku$ with
$\BP\<1\>$ as a wedge summand. In~\cite[1.18]{HKM} it is conjectured that 
the core of $\ku$ is $\BP\<1\>$. Here $\ku=\BP\<1\>$ if $p=2$. Since 
$H^*(\BP\<1\>;\bF_p)$ is a cyclic $A$-module, this is now immediate.

\begin{exmp}\label{prop:BP<1>-Core}
The spectrum $\BP\<1\>$ is minimal atomic, hence the canonical mono\-mor\-phism 
$\BP\<1\>\rtarr\ku$ is a core.
\end{exmp}

More generally, $H^*(\BP\<n\>)$ is a cyclic $A$-module for all $n\geqslant -1$, the 
extreme cases being $\BP\<-1\>=H\F_p$ and $\BP\<0\>=H\Z_{(p)}$.

\begin{exmp}\label{thm:Core-BP<m>}
For $n\geqslant -1$, $\BP\<n\>$ is a minimal atomic spectrum.
\end{exmp}

The following example of the non-uniqueness of cores generalizes \cite[1.17]{HKM}.

\begin{exmp}\label{cor:Cores-different} For $n\geqslant 0$, the canonical maps 
$$
\xymatrix@1{
\BP\ar[r] &\BP\wedge\BP\<n\>&\BP\<n\>\ar[l]}$$
induced by the units of $\BP$ and $\BP\<n\>$ are both cores of $\BP$.
\end{exmp}
\begin{proof}
The left map is a monomorphism since it factors the ($p$-local)
Hurewicz homomorphism of $\BP$. The right map is a monomorphism
since it is split by the $\BP$-action $\BP\sma \BP\<n\>\rtarr \BP\<n\>$.
\end{proof} 

Since $H^*(\ko;\F_2) = A//A(1)$ and $H^*(\eo_2;\F_2) = A//A(2)$ are cyclic 
$A$-modules, we have the following complement to Proposition \ref{prop:BP<1>-Core}.

\begin{prop} At $p=2$, $\ko$ and $\eo_2$ are minimal atomic spectra.
\end{prop}

Some well-known Thom complexes give further examples.

\begin{prop} Let $X$ be $\RP^\infty_{-1}$, $\CP^\infty_{-1}$, 
or $\HP^\infty_{-1}$, that is, the Thom spectrum of the negative 
of the canonical real, complex, or quaternionic line bundle.
At $p=2$, $X$ is minimal atomic.
\end{prop}
\begin{proof}
Let $d=1$, $2$, and $4$ and $P=\RP^\infty$, $\CP^\infty$, and
$\HP^\infty$ in the respective cases. Then $H^*(P;\F_2)=\F_2[x]$, 
where $x\in H^d(P;\F_2)$ is the $d$th Stiefel-Whitney class of the
canonical line bundle. Since $X$ is a Thom spectrum, $H^*(X;\bF_2)$ 
is the free $H^*(P;\F_2)$-module generated by the Thom class $\mu$
in degree $-d$. A standard calculation shows that $\mathrm{S}q^{nd}\mu=x^n\mu$ 
for $n\geqslant 1$, so $H^*(X;\F_2)$ is cyclic over $A$.
\end{proof}

To give some examples where we must check the ``no mod $p$\, homotopy'' 
condition directly, we consider a few suspension spectra and another Thom 
spectrum. Let $\xi_3\rtarr \HP^\infty$ be the bundle associated to the 
adjoint representation of $S^3$ and let $M\xi_3$ be its Thom complex 
(also known as a quaternionic quasi-projective space). It has one cell in 
each positive dimension congruent to $3\pmod{4}$.

By~\cite{Holzager, Cooke&Smith,7A}, for each odd prime~$p$, there 
is a splitting of $p$-local spaces
$$\Sigma \CP^\infty\simeq
W_{1}\vee W_{2}\vee\cdots\vee W_{p-1},$$
where $W_{r}$ has cells in all dimensions of the form
$2(p-1)k+2r+1$ with $k\geqslant 0$.

\begin{prop}\label{prop:Priddy-NuclearTest-Ex}
At the prime~$2$, $\Sigma^\infty\CP^\infty$, $\Sigma^\infty\HP^\infty$
and $\Sigma^\infty M\xi_3$ are minimal atomic spectra. At an odd prime~$p$, 
each $\Sigma^\infty W_{r}$ is minimal atomic.
\end{prop}
\begin{proof}
Let $a(n)=1$ if $n$ is even and $a(n)=2$ if $n$ is odd.
By ~\cite{DSegal}, the Hurewicz homomorphisms
$$ h\:\pi_{2n}(\Sigma^\infty\CP^\infty)\rtarr H_{2n}(\CP^\infty)\iso\Z
\ \ \ \text{and}\ \ \
h\:\pi_{4n}(\Sigma^\infty\HP^\infty)\rtarr H_{4n}(\HP^\infty)\iso\Z$$
have images of index $n!$ and $(2n)!/a(n)$, respectively. Thus, for $n>1$, 
the corresponding mod~$2$ Hurewicz homomorphisms are trivial. By~\cite{Walker},
the Hurewicz homomorphism
$$h\:\pi_{4n+3}(\Sigma^\infty M\xi_3)\rtarr H_{4n+3}(M\xi_3)\iso\Z$$
has image of index $a(n)(2n-1)!$ , so for each $n\geqslant1$ the
associated mod~$2$ Hurewicz homomorphism is also trivial. 
The odd primary results follow similarly from the calculation of $h$ for 
$\Sigma^\infty\CP^\infty$.
\end{proof}

\begin{rem}\label{rem:RP}{\ } We raise a few questions here.
\begin{enumerate}
\item[(i)] There are many basic results in the literature in which
interesting spaces are split $p$-locally into products of indecomposable 
factors and interesting spectra are split $p$-locally into wedges of indecomposable
summands. (The notion of wedge indecomposability is less interesting in the case of
spaces). It is a very interesting set of problems to revisit these splittings and 
determine which of the summands are atomic rather than just indecomposable, and
which are minimal atomic rather than just atomic. The results above just give 
particularly elementary examples. 
\item[(ii)]  The suspension spectrum of $\RP^\infty$ presents an 
interesting challenge. It is a standard observation that 
$H^*(\RP^\infty;\F_2)$ is an atomic, but not cyclic, $A$-module, 
in the sense that any $A$-endomorphism which is the identity on 
$H^1(\RP^\infty;\F_2)$ is an isomorphism. This implies that 
$\Sigma^\infty\RP^\infty$ is atomic. However, since the top cell
of $\RP^3$ splits off stably, the stable Hurewicz homomorphism 
$\pi_3(\Sigma^\infty\RP^\infty)\rtarr H_3(\RP^\infty;\F_2)$ is 
non-trivial, hence $\Sigma^\infty\RP^\infty$ cannot be minimal 
atomic. It would be interesting to identify a core of 
$\Sigma^\infty\RP^\infty$. For an odd prime~$p$, similar remarks 
apply to the $(p-1)$ wedge summands of $\SI^{\infty}B\Z/p$, one of 
which is $(B\Sigma_p)_{(p)}$.
\item[(iii)] Fred Cohen reminds us that it is an open question whether or not $K(\bZ/2,n)$
or the $p-1$ summands of $\Sigma K(\bZ/p,n)$ are stably atomic for $n\geqslant 2$. The
question was posed by Priddy in \cite[p.\, 379]{Cont}.
\end{enumerate}
\end{rem}

\section{A construction of the spectrum $\BoP$}\label{sec:BoP}

In this section, all spectra are understood to be localized at~$2$, 
and $S=S^0$. Recall the spectrum $\BoP$ of Pengelley~\cite{Pengelley}. 
It has no mod $2$\, detectable homotopy \cite[5.5]{Pengelley} and 
it is a retract of $\MSU$, so we have a monomorphism $j: \BoP\rtarr \MSU$. 

\begin{exmp}\label{prop:BoP}
The monomorphism $j\:\BoP\rtarr\MSU$ is a core of $MSU$.
\end{exmp}

We recall a further property of $\BoP$, proven in Pengelley \cite[6.15, 6.16]{Pengelley}.

\begin{prop}\label{Pen} 
There is a map $p: \BoP\rtarr \ko$ that induces an
epimorphism on homotopy groups
in all degrees and an isomorphism in odd degrees.
\end{prop}

\begin{cor} The odd degree homotopy groups of the fiber $Fp$ are zero.
\end{cor}

We now give a description of $\BoP$ as a nuclear spectrum, thus providing a 
simple construction of it that is independent of \cite{Pengelley}. Guided
by Proposition \ref{Pen}, we construct a nuclear spectrum $X$ and a map
$q: X\rtarr \ko$ that induces a monomorphism on homotopy groups in odd
degrees, and we prove that it induces an epimorphism on homotopy groups.
That turns out to imply that $X$ is equivalent to $\BoP$.

We begin with $X_0=S$, and we inductively attach even dimensional cells, 
letting $X_{2n} = X_{2n+1}$ for all $n\geqslant 0$. Suppose that we have factored 
the unit $\io: S\rtarr\ko$ through a map $q_n\: X_{2n-1}\rtarr\ko$. We enlarge 
$X_{2n-1}$ to $X_{2n}$ by attaching $2n$-cells minimally, so that (\ref{key}) is 
satisfied. We do this so as to kill the kernel of 
$${q_n}_*: \pi_{2n-1}(X_{2n-1})\rtarr\pi_{2n-1}(\ko).$$
Thus, in the resulting cofiber sequence
$$ J_{2n-1}\rtarr X_{2n-1}\rtarr X_{2n},$$
$$\text{Im}(\pi_{2n-1}(J_{2n-1})\rtarr\pi_{2n-1}(X_{2n-1}))=
\text{Ker}(\pi_{2n-1}(X_{2n-1})\rtarr\pi_{2n-1}(\ko)).$$
Clearly $q_n$ extends to a map
$$q_{n+1}\: X_{2n} = X_{2n+1}\rtarr\ko.$$ 
In the limit we obtain a nuclear complex $X$ and a map $q\:X\rtarr\ko$ that induces an 
isomorphism on $\pi_0$ and a monomorphism on $\pi_*$ in odd degrees.

\begin{prop}\label{prop:Fibre(BoP'->\ko)}
$q\: X\rtarr\ko$ induces an epimorphism on homotopy groups.
\end{prop}

\begin{cor} The odd degree homotopy groups of the fiber $Fq$ are zero.
\end{cor}

Let $\nu\in\pi_3(S)$ and $\sigma\in\pi_7(S)$ be the Hopf maps. If $x\in X$ 
has even degree, then $\nu x$ and $\si x$ are odd degree elements of the kernel 
of $q_*$, hence they are zero. The proposition is therefore a direct consequence
of the following result, which is presumably known. Since we do not know of 
a reference for it, we will give a proof in the next section.

\begin{prop}\label{prop:S->X->\ko}
Let $X$ be a Hurewicz complex of dimension $0$ with inclusion
of the bottom cell $i\:S\rtarr X$ and let $q\:X\rtarr\ko$ be a map
such that the composite $S\xrightarrow{i}X\xrightarrow{q}\ko$ is 
the unit $\io: S\rtarr \ko$. If $\nu x=0$ and $\sigma x=0$ in 
$\pi_*(X)$ for every even degree element $x\in\pi_*(X)$, then 
$q_*\:\pi_*(X)\rtarr\pi_*(\ko)$ is an epimorphism.
\end{prop}

\begin{thm}\label{thm:BoP'=BoP}
There are equivalences $f\:X\rtarr \BoP$ and $g\:\BoP \rtarr X$ such that
the following diagram is homotopy commutative. 
$$\xymatrix{
X\ar[r]^-{f}\ar[dr]_{q}& \BoP\ar[r]^-{g}\ar[d]^{p}& X\ar[dl]^{q} \\
&\ko&}$$
\end{thm}
\begin{proof}
We construct maps $f$ and $g$ such that the diagram is homotopy commutative.
The maps $f$ and $g$, hence also the composites $g\com f$ and $f\com g$, then 
induce isomorphisms on $\pi_0$. Since $X$ and $\BoP$ are atomic, these composites 
are equivalences and therefore $f$ and $g$ are equivalences. We may take $\BoP$ 
and $\ko$ to be Hurewicz complexes and take $p$ to be the identity map on the 
bottom cell. Taking $f_0\:X_0=X_1=S\rtarr \BoP$ to be the identity map on the bottom cell
and $h_0$ to be the constant homotopy at the identity map, we assume inductively that we 
have a map $f_n\:X_{2n-1}\rtarr \BoP$ and a homotopy $h_n\:q_n\htp p\com f_n$. 
Consider the following diagram, where we implicitly
precompose maps already specified with the map of cells $CJ_{2n-1}\rtarr X_{2n+1}$
that constructs $X_{2n+1}$ from $X_{2n-1}$.
$$\xymatrix{
J_{2n-1} \ar[dd] \ar[rr]^-{i_0} & & J_{2n-1} \sma I_+ \ar[dl]_-{h_{n}} \ar[dd] 
& & J_{2n-1} \ar[ll]_-{i_1} \dlto_-{f_n} \ar[dd] \\
& \ko & & \BoP \ar[ll]_(.3){p} & \\
CJ_{2n-1} \ar[rr]_-{i_0} \ar[ur]^-{q_{n+1}} 
&& CJ_{2n-1}\sma I_+ \ar@{-->}[ul]^-{h_{n+1}} & & 
CJ_{2n-1} \ar@{-->}[ul]^-{f_{n+1}} \ar[ll]^-{i_1} \\}$$
Since $J_{2n-1}$ is a wedge of $(2n-1)$-spheres and $\pi_{2n-1}(Fp) = 0$,
$[J_{2n-1},Fp] = 0$. A standard result, given in just this form in \cite[Lemma 1]{May}, 
shows that there are maps $f_{n+1}$ and $h_{n+1}$ that make the diagram commute.
Passing to colimits, we obtain $f$ and a homotopy $h\:q\htp p\com f_n$. Since the
homology groups of $\BoP$ are concentrated in even degrees \cite{Pengelley}, we 
can replace it by a minimal complex, with cells only in even degrees. This allows
us to reverse the roles of $X$ and $\BoP$ to construct $g$.
\end{proof}

A similar argument proves the following result.

\begin{prop}\label{prop:MSU->BoP'}
There is a map $r\:\MSU\rtarr X$ such that the following diagram is homotopy commutative.
$$\xymatrix{
\MSU\ar[rr]^{r}\ar[dr]_{t}&&X\ar[dl]^{q} \\
&\ko &}$$
\end{prop}

It is not clear that $\BoP$ is the only core of $\MSU$ up to equivalence, 
but we conjecture that it is. The following consequence of 
Lemma~\ref{prop:S->X->\ko} may shed some light on this question.

\begin{prop}\label{prop:MSU-Core->\ko}
If $Y\rtarr\MSU$ is a core, the composite $Y\rtarr\MSU\rtarr\ko$ 
induces an epimorphism on homotopy groups.
\end{prop}

\begin{rem} It might be of interest to revisit the results of~\cite{Kochman-BoP, Pengelley}
from our present perspective. However, it is not clear how to construct a map 
$X\rtarr\MSU$ that induces the identity on $\pi_0$ and how
the distinguished map of~\cite{Kochman-BoP} fits in. It might be of more interest
to revisit the results of \cite {Kochman-BoP, Pengelley} from the perspective of
$S$-modules \cite{EKMM}. Pengelley constructs $\BoP$ by first constructing 
another spectrum, which he denotes by $X$, and then taking a fiber to kill $\BP$ summands in it. 
His $X$ is obtained from $\MSU$ by using the Baas-Sullivan theory of manifolds with singularities 
to kill a regular sequence of elements in $\pi_*(MSU)$. We can instead use the results of 
\cite[Ch.\,V]{EKMM} to construct $X$ as an $\MSU$-module together with a map of $\MSU$-modules 
$\MSU\rtarr X$. It seems plausible that the methods of \cite{EKMM, Strick} can be used to construct $\BoP$
as a commutative $\MSU$-ring spectrum.
\end{rem}

\section{The proof of Proposition \ref{prop:S->X->\ko}}\label{sec:\ko}

We continue to work with spectra localized at $2$. Recall that
\begin{equation}\label{eqn:\ko_*}
\pi_*(\ko)=\Z_{(2)}[\et,\al,\be]/(2\et,\et^3,\et\al,\al^2-4\be),
\end{equation}
where $\deg \et=1$, $\deg \al=4$, and $\deg \be =8$. We will describe elements 
of $\pi_*(X)$ that map to each of the additive generators of $\pi_*(\ko)$. Note that,
since we do not know that $X$ is a ring spectrum, we cannot exploit the algebra
structure of $\pi_*(\ko)$. The essential point is to describe additive generators in 
terms of Toda brackets in $\pi_*(\ko)$ that admit analogues in $\pi_*(X)$. 

We are interested in Toda brackets of the form $\<a,b,c\>$, where $a$ 
and $b$ are elements of $\pi_*(S)$ and $c$ is an element of $\pi_*(Y)$ for a 
spectrum $Y$. We require $ab=0$ and $bc=0$, and then $\<a,b,c\>$ is a coset
of elements in $\pi_{|a|+|b|+|c|+1}(Y)$ with respect to the indeterminacy subgroup
\begin{align*}
\indeter\<a,b,c\>&=a\pi_{|b|+|c|+1}(Y)+(\pi_{|a|+|b|+1}(S))c.
\end{align*}
Such Toda brackets are natural with respect to maps $Y\rtarr Z$.

\begin{rem} We remark parenthetically that the theory of Toda brackets simplifies
greatly if one defines them in terms of the associative smash product in one of the
modern categories of spectra, such as the category of $S$-modules of \cite{EKMM}. 
A systematic exposition would be of value. In brief, the conclusion must be that 
all of the results that are catalogued in \cite{May0} for matric Massey products 
in the homology of DGA's carry over verbatim to $S$-modules. 
\end{rem}

Now take $X$ as in Proposition \ref{prop:S->X->\ko}. Recall that $8\nu=0$ and $16\si=0$ 
in $\pi_*(S)$ and that, by hypothesis, $\nu$ and $\si$ annihilate all even degree 
elements of $\pi_*(X)$. Let $b_0$ denote $i:S\rtarr X$ regarded as an element of $\pi_0(X)$
and choose coset representatives in iterated Toda products as follows:
$$a_1\in\<8,\nu,b_0\>,\ \ \ b_{k}\in\<16,\sigma,b_{k-1}\>,
\ \ \text{and}\ \ \ a_{k+1}\in\<16,\sigma,a_{k}\>,$$
where $k\geqslant 1$. The indeterminacies are benign for our purposes since they are
\begin{align*}
\indeter a_1&=(\pi_4(S))b_0+8\pi_4X=8\pi_4(X), \\
\indeter b_k&=(\pi_8(S))b_{k-1}+16\pi_{8k}(X)\equiv 16\pi_{8k}(X)\ \text{mod}\Ker(q_*) \\
\indeter a_k&=(\pi_8(S))a_{k-1}+16\pi_{8k-4}(X)\equiv 16\pi_{8k-4}(X) \ \text{mod}\Ker(q_*).
\end{align*}
Here the congruences hold since $\pi_8(S)$ is $2$-torsion and there are no torsion 
elements in the relevant degrees of $\pi_*(\ko)$. For $k\geqslant 0$, we also have the elements
$$\mu_{8k+1}b_0\in\pi_{8k+1}(X) \ \ \text{and}\ \ \mu_{8k+2}b_0\in\pi_{8k+2}(X),$$
where $\mu_{8k+1}$ and $\mu_{8k+2}$ are the usual elements in $\pi_*(S)$. 
Now $q_*\:\pi_*(X)\rtarr\pi_*(\ko)$ maps these elements to elements
of the same form in $\pi_*(\ko)$, where $b_0\in\pi_0(\ko)$ is the unit of $\ko$.
In the familiar periodic pattern $\bZ_2$, $\bZ_2$, $0$, $\bZ$, $0$, $0$, $0$, $\bZ$, 
the additive positive degree generators of $\pi_*(\ko)$ are 
$$\et\be^k = \mu_{8k+1}b_0,\ \ \et^2\be^k = \mu_{8k+2}b_0,\ \ \al\be^{k},\ \ \text{and}\ \ \be^{k+1},$$ 
where $k\geqslant 0$. The following known result gives that $\al\be^k = a_{k+1}$ and 
$\be^{k+1}= b_{k+1}$ in $\pi_*(\ko)$, and this completes the proof that $q_*$ is an epimorphism.

\begin{lem} In $\pi_*(\ko)$,
$$\al \in\<8,\nu,b_0\>, \ \  \be^{k}\in\<16,\sigma,\be^{k-1}\>,
\ \ \text{and}\ \  \al \be^{k}\in\<16,\sigma,\al\be^{k-1}\>$$
for $k\geqslant 1$, where the indeterminancy is $0\,\text{mod}\, 2$ in each case.
\end{lem}

An unstable version of the lemma is stated without proof in ~\cite[p.64]{Whitehead}, 
where it is attributed to Barratt. One quick way to see the result is to use the convergence
of Massey products to Massey products in the May spectral sequence and of Massey products to 
Toda brackets in the Adams spectral sequence, but the details would take us too far afield.

\section*{Appendix: Errata to \cite{HKM}}

We take this opportunity to correct some minor errors in the 
proof of ~\cite[2.11]{HKM}. In brief, the last 
two sentences of the cited proof should be replaced with the 
following two sentences.
``If $p=2$, then $Q^8(a_1)\equiv a_5$ mod decomposables, and, if $p>2$, then 
$Q^{2p}(a_{p-1})\equiv a_{(2p+1)(p-1)}$ mod decomposables, by \cite{P0} or 
\cite[II.8.1]{CLM}. Here $a_{p-1}$ is in the image of $H_*(BP)$, but $H_*(BP)$ 
has no indecomposable elements in degree $10$ if $p=2$ or in
degree $2(2p+1)(p-1)$ if $p>2$.''

\end{document}